\documentclass[12pt,reqno]{amsart}
\usepackage{amssymb}

\theoremstyle{plain}
\newtheorem{theorem}{Theorem}
\newtheorem{corollary}[theorem]{Corollary}
\newtheorem{lemma}[theorem]{Lemma}

\newtheorem*{problem}{Problem}

\theoremstyle{definition}

\theoremstyle{remark}

\def\Bl1{\mathcal B(\ell_1)}
\def\lin{{\rm lin}}
\def\mid{\::\:}

\newcommand{\marg}[1]{}
\newcommand{\lab}[1]{\label{#1}}
\newcommand{\abs}[1]{\lvert#1\rvert}
\newcommand{\norm}[1]{\lVert#1\rVert}
\renewcommand{\le}{\leqslant}
\renewcommand{\ge}{\geqslant}

\begin{document}
\baselineskip 18pt

\title{On the modulus of C.~J.~Read's operator}

\author[V.~G.~Troitsky]{Vladimir~G.~Troitsky}
\address{
         Mathematics Department\\
         University of Illinois at Urbana-Champaign\\
         1409 West Green St\\
         Urbana, IL 61801, USA}
\email{vladimir@math.uiuc.edu}

\thanks{Supported in part by NSF Grant DMS 96-22454.}

\date{February 13, 1998}

\begin{abstract}
  Let $T\colon\ell_1\to\ell_1$ be the quasinilpotent operator without
  an invariant subspace constructed by C.~J.~Read in~\cite{R3}. We
  prove that the modulus of this operator has an invariant subspace
  (and even an eigenvector). This answers a question posed by
  Y.~Abramovich, C.~Aliprantis and O.~Burkinshaw in~\cite{AAB1,AAB3}.
\end{abstract}

\maketitle

During the last several years there has been a noticeable increase of
interest in the invariant subspace problem for positive operators on
Banach lattices. A rather complete and comprehensive survey on this
topic is presented in~\cite{AAB3}, to which we refer the reader for
details and for an extensive bibliography.  In particular, the
following theorem was proved in~\cite{AAB1}.
\begin{theorem}[\cite{AAB1,AAB3}]  \lab{t:aab}
  If the modulus of a continuous operator $T\colon\ell_p\to\ell_p$
  $(1\le p<\infty)$ exists and is quasinilpotent, then $T$ has a
  non-trivial closed invariant subspace which is an ideal.
\end{theorem}
It follows that each positive quasinilpotent operator on $\ell_p \ 
(1\le p<+\infty)$ has a nontrivial closed invariant subspace. In the
same papers the authors posed the following problem.
\begin{problem}
  Does every positive operator on $\ell_1$ have an invariant subspace?
\end{problem}

Keeping in mind that each operator on $\ell_1$ has a modulus and
that C.~J.~Read in~\cite{R1,R2,R3} has constructed several operators
on $\ell_1$ without invariant subspaces, it was suggested
in~\cite{AAB1,AAB3} that the modulus of some of these operators might
be a natural candidate for a counterexample to the above problem.
Following this suggestion, we will be dealing in this paper with the
modulus of the quasinilpotent operator $T$ constructed in ~\cite{R3}.
It turns out, quite surprisingly, that not only does $\abs{T}$ have an
invariant subspace but it even has a positive eigenvector. This result
increases the chances for an affirmative answer to the above problem.

The paper is organized as follows.  After introducing some
necessary notation and terminology we prove a general theorem on the
existence of an invariant subspace for the modulus of a
quasinilpotent operator.  The rest of the paper will be devoted to the
verification that C.~J.~Read's operator, constructed in~\cite{R3},
satisfies all the hypotheses of this theorem and so its modulus does
have an invariant subspace.

For terminology and notation regarding operators and Banach lattices we
refer to~\cite{AB,Sch}.  All operators considered in this work are linear
and continuous. The space of all operators on a Banach space $X$ is
denoted by ${\mathcal B}(X)$, while ${\mathcal K}(X)$ stands for the
subspace of all compact operators. A linear operator on a Banach
lattice is said to be \textit{positive} if it maps positive vectors to
positive vectors. By an \textit{invariant subspace} (\textit{invariant
  ideal}) of an operator we mean a closed nontrivial subspace (resp.
closed nontrivial ideal) which is invariant under the operator.

Together with the usual operator norm $\norm{S}$ of an operator
$S\colon X\to X$
on a Banach space, we will also consider
the (essential) seminorm $\norm{S}_e$ given by
$$
  \norm{S}_e=\inf\{\,\norm{S-K}\mid K\in {\mathcal K}(X)\,\}.
$$
The \textit{essential spectral radius} of $S$ is computed via the
formula $r_e(S)=\lim_n\sqrt[n]{\norm{S^n}_e}$. This, of course, is an
analogue of the familiar formula for the usual
spectral radius $r(S)=\lim_n\sqrt[n]{\norm{S^n}}$.

It is obvious that $\norm{S}_e\le\norm{S}$ and $r_e(S)\le r(S)$. It is
known that if $r_e(S)=0$, then every nonzero point of $\sigma(S)$ is an
eigenvalue. Further details on essential spectral radius can be found
in~\cite{N1,CPY}.  We will use the following important version of the
Krein-Rutman theorem established by R.~Nussbaum.

\begin{theorem}{\cite[Corollary~2.2]{N2}} \lab{t:kr}
  Let $S$ be a positive operator on a Banach lattice such that
  $r_e(S)<r(S)$, then $r(S)$ is an eigenvalue of $S$ corresponding to
  a positive eigenvector.\footnote{In~\cite{N2} a more general form of
  this theorem is given which is valid for ordered Banach spaces.}
\end{theorem}

We use this fact in the proof of the following simple but
rather unexpected result.

\begin{theorem} \lab{t:general}
  Suppose that a quasinilpotent operator $S$ on $\ell_p$ has no
  invariant ideals and $S^-$ is compact.  Then $r(\abs{S})$
  is a positive eigenvalue of $\abs{S}$ corresponding to a positive
  eigenvector. In particular, $\abs{S}$ has an invariant
  subspace.
\end{theorem}

\begin{proof}
  First observe that the operator $\abs{S}$ cannot be quasinilpotent.
  Indeed, if it were, then by Theorem~\ref{t:aab}
   the operator $S$ itself would have an
  invariant closed ideal contrary to our hypothesis. Thus,
  $r(\abs{S})>0$.

  Next we claim that $r_e(\abs{S})=0$. To prove this, notice that
  $\abs{S}=S+2S^-$, and so
  \[
    \abs{S}^n=(S+2S^-)^n=S^n+RS^-,
  \]
  where $R$ is some polynomial in $S$ and $S^-$. Hence $RS^-$ is compact,
  whence
  \[
    \norm{\abs{S}^n}_e=\norm{S^n+RS^-}_e\le\norm{S^n},
  \]
  and consequently
  \[
    r_e(\abs{S})\le\lim_{n\to\infty}\sqrt[n]{\norm{S^n}}=r(S)=0.
  \]
An application of Theorem~\ref{t:kr} finishes the proof.
\end{proof}

\begin{corollary}
Under the hypotheses of the above theorem the operator $S^+$
also has a nontrivial closed invariant subspace.
\end{corollary}
\begin{proof}
  There are two possibilities: either $S^+$ is quasinilpotent or it is
  not.  If $S^+$ is quasinilpotent, then applying Theorem~\ref{t:aab}
  again, we see that $S^+$ has an invariant ideal.
  
  Assume that $S^+$ is not quasinilpotent. Since $S^+ =S + S^-$, the
  same argument as in the proof of Theorem~\ref{t:general} shows that $r_e(S^+)=0$
  and we can again apply Theorem~\ref{t:kr}.
\end{proof}

Theorem~\ref{t:general} is strong enough to enable us to prove
Corollary~\ref{c:main} about the modulus of C.~J.~Read's operator.
But first we would like to mention a nice generalization of
Theorem~\ref{t:general}. Recall that a positive operator
on a Banach lattice is called
\textit{compact-friendly} if it commutes with another positive
operator which dominates some non-zero operator which in turn is
dominated by a positive compact operator. The class of
compact-friendly operators was introduced and studied by Abramovich,
Aliprantis, and Burkinshaw in~\cite{AAB2}.
This class includes positive operators that dominate or are dominated
by a non-zero compact operator, or commute with a non-zero compact
operator. Also, every positive operator on any discrete Banach lattice
and every positive kernel operator is compact-friendly. It
follows from \cite[Theorem 11.2]{AAB3} that if a positive
compact-friendly operator is quasinilpotent, then it has an invariant
ideal. Mimicking the proof of Theorem~\ref{t:general} we can obtain
the following theorem.

\begin{theorem}
Let $B$ be a  compact-friendly operator  with $r_e(B)=0$. Then $B$ has a
nontrivial invariant subspace.
\end{theorem}

To illustrate this theorem we mention the following result:
{\it If $S$ is a quasinilpotent kernel operator  and $S^-$ (resp. $S^+$) is
compact, then  $\abs{S}$ and $S^+$ (resp. $S^-$) have invariant
subspaces}.

\bigskip

Recall that an operator $S\colon X\to Y$ between two Banach spaces is
called \textit{nuclear} if it can be written in the form
$S=\sum_{i=0}^\infty x_i^*\otimes y_i$ with $x_i^*$ in $X^*$, $y_i$ in
$Y$ and $\sum_{i=0}^\infty \norm{x_i^*}\norm{y_i}<\infty$. Here, as
usual, the elementary tensor $x^*\otimes y: X\to Y$ is given by
$(x^*\otimes y)(x)=x^*(x)y$. The \textit{nuclear norm} $\nu(S)$ is
defined by $\nu(S)=\inf\sum_{i=0}^\infty\norm{x_i^*}\norm{y_i}$, where
the infimum is taken over all nuclear representations of $S$. 
For a nuclear operator $S$ we have
$\norm{S}\le\nu(S)$. This implies in particular that every nuclear
operator can be approximated by finite-rank operators and, therefore,
is compact.

Following \cite{R1,R2,R3} we denote the standard unit vectors of
$\ell_1$ by $(f_i)_{i=0}^\infty$.  It is well known that we can
consider each $S\in\Bl1$ as an infinite matrix
$S=(s_{ij})_{i,j=0}^\infty$. Let $S_{(i)}$ denote the $i$-th row of
this matrix. If $x\in\ell_1$, then $(Sx)_i=\langle S_{(i)},x\rangle$,
so that $Sx=\sum_{i=0}^\infty\langle S_{(i)},x\rangle f_i$. This gives
a nuclear representation $S=\sum_{i=0}^\infty S_{(i)}\otimes f_i$,
where the rows $S_{(i)}$ of $S$ are considered as linear functionals
on $\ell_1$. It follows that
\[
  \nu(S)\le\sum_{i=0}^\infty \norm{S_{(i)}}_\infty\norm{f_i}_1=
  \sum_{i=0}^\infty \norm{S_{(i)}}_\infty,
\]
so that $S$ is nuclear if the last sum is finite.

In spite of the fact that the construction of a quasinilpotent operator
on $\ell_1$ without an invariant subspace in~\cite{R3} is far
from being simple, the presentation in~\cite{R3} is very clearly
structured.  In Sections 2 and 3 of~\cite{R3} C.~J.~Read presents the
construction of operator $T$, and in the rest of the paper he proves
that $T$ is bounded (Lemma~5.1), quasinilpotent (Theorem~6.5), and has
no invariant subspaces (Theorem~7). In what follows we will restate
(practically verbatim and retaining the notation) the definition of
$T$. In the proof of~\cite[Lemma~5.1]{R3} C.~J.~Read reveals a lot of
information about the structure of the infinite matrix of $T$. Since
we also are interested in this structure, we incorporate some
fragments of the proof of~\cite[Lemma~5.1]{R3} in our proof
of~Lemma~\ref{l:nuclear}.

The symbol $F$ denotes the linear subspace of $\ell_1$, spanned by
$f_i$'s, and thus, $F$ is dense and consists of all eventually vanishing
sequences. Let ${\bf d}=(a_1,b_1,a_2,b_2,\ldots)$ be a strictly increasing
sequence of positive integers. Also let $a_0=1$, $v_0=0$,
and $v_n=n(a_n+b_n)$ for $n\ge 1$. Then there is a unique sequence
$(e_i)_{i=0}^\infty\subset F$ with the following properties:
\begin{enumerate}
  \item[0)] $f_0=e_0$;
  \item[A)] if integers $r$, $n$, and $i$ satisfy $0<r\le n$,
    $i\in[0,v_{n-r}]+ra_n$, then
    $f_i=(n^{ra_n}e_i-e_{i-ra_n})(n-r)^{i-ra_n}a_{n-r}$;
  \item[B)] if integers $r$, $n$, and $i$ satisfy $0<r<n$,
    $i\in(ra_n+v_{n-r},(r+1)a_n)$, (respectively, $1\le n$,
    $i\in(v_{n-1},a_n))$, then $f_i=n^i2^{(h-i)/\sqrt{a_n}}e_i$, where
    $h=(r+\frac{1}{2})a_n$ (respectively, $h=\frac{1}{2}a_n$);
  \item[C)] if integers $r$, $n$, and $i$ satisfy $0<r\le n$,
    $i\in[r(a_n+b_n),na_n+rb_n]$, then
    $f_i=n^ie_i-b_nn^{i-b_n}e_{i-b_n}$;
  \item[D)] if integers $r$, $n$, and $i$ satisfy $0\le r<n$,
  $i\in(na_n+rb_n,(r+1)(a_n+b_n))$, then $f_i=n^i2^{(h-i)/\sqrt{b_n}}e_i$,
  where $h=(r+\frac{1}{2})b_n$.
\end{enumerate}

Indeed, since $f_i=\sum_{j=0}^{i}\lambda_{ij}e_j$ for each $i\ge 0$
and $\lambda_{ii}$ is always nonzero, this linear relation is
invertible. Further,
\begin{equation*}
  \lin\{e_i\mid i=1,\dots,n\}=
  \lin\{f_i\mid i=1,\dots,n\}\mbox{ for every }n\ge 0.
\end{equation*}
In particular all $e_i$ are linearly
independent and span $F$. Then C.~J.~Read defines $T\colon F\to F$ to be
the unique linear map such that $Te_i=e_{i+1}$, and in Lemma~5.1
he proves that $\norm{Tf_i}\le1$ for every $i\ge 0$ provided {\bf d} increases
sufficiently rapidly, i.~e.,  satisfies several conditions of the form
\begin{eqnarray*}
  a_n & \ge & G(n, a_0, b_0, a_1, b_1, \dots, a_{n-1}, b_{n-1}),
  \mbox{ and} \\
  b_n & \ge & H(n, a_0, b_0, a_1, b_1, \dots, a_{n-1}, b_{n-1}, a_n),
\end{eqnarray*}
where $G$ and $H$ are some positive integer-valued
functions.  It follows that $T$ can be extended to
a bounded operator on $\ell_1$. Finally, in Theorems~6.5 and~7, C.~J.~Read
proves that this extension, which is also denoted by $T$, is
quasinilpotent and has no invariant subspaces, provided {\bf d}
increases sufficiently rapidly.

Our plan is as follows: we will prove that the negative part of $T$ is
nuclear, hence compact. Then Theorem~\ref{t:general} will imply that
$r(\abs{T})$ is a positive eigenvalue of $\abs{T}$, corresponding to a
positive eigenvector.

\begin{lemma} \lab{l:nuclear}
  The operator $T^-$ is nuclear, provided {\bf d} increases
  sufficiently rapidly.
\end{lemma}

\begin{proof}
  Similarly to the proofs of~\cite[Lemma~5.1]{R3}
  and~\cite[Lemma~6.1]{R2} we study the matrices
  $(t_{ij})_{i,j=0}^\infty$ and $(t^-_{ij})_{i,j=0}^\infty$ of $T$ and
  $T^-$ respectively.
  Recall that $t_{ki}=(Tf_i)_k$, so that it suffices to look at the
  images of the standard unit vectors under $T$. We will see that the
  matrix of $T$ is quite sparse and has the following structure: every
  entry on the diagonal right under the main diagonal is strictly
  positive, there is no nonnegative entries below this diagonal, and
  there are some entries above it. We consider consecutively all the
  cases mentioned above.

\begin{enumerate}
  \item[0)] $Tf_0=e_1=2^{(1-a_1/2)/\sqrt{a_1}}f_1$, so that $T^-f_0=0$.
  \item[A)] If $i<v_{n-r}+ra_n$, i.~e. $i$ is not the right end point
    of the interval $[ra_n, v_{n-r}+ra_n]$, then $Tf_i=(n-r)^{-1}f_{i+1}$, so
    that $T^-f_i=0$. The only nontrivial case here is when $i$ is the
    right end of the interval, i.~e. $i=v_{n-r}+ra_n$. Then we have
    \begin{multline*}
      Tf_i=
      a_{n-r}n^{ra_n}(n-r)^{v_{n-r}}e_{1+ra_n+v_{n-r}}-
      a_{n-r}(n-r)^{v_{n-r}}e_{1+v_{n-r}}\\
      =\varepsilon_1f_{1+v_{n-r}+ra_n}-\varepsilon_2f_{1+v_{n-r}},
    \end{multline*}
    where $\varepsilon_1>0$ and $\varepsilon_2$ is given by
    \[
      \varepsilon_2=(n-r+1)^{-1-v_{n-r}}
        2^{(1+v_{n-r}-a_{n-r+1}/2)/\sqrt{a_{n-r+1}}}
        a_{n-r}(n-r)^{v_{n-r}},
    \]
    so that
    \marg{e:A}
    \begin{equation} \label{e:A}
      T^-f_{v_{n-r}+ra_n}=\varepsilon_2f_{1+v_{n-r}}.
    \end{equation}
  \item[B)] Similarly, if $ra_n+v_{n-r}<i<(r+1)a_n-1$ or $v_{n-1}<i<a_n-1$,
    then
    $Tf_i=n^{-1}2^{1/\sqrt{a_n}}f_{i+1}$, so that $T^-f_i=0$.
    If $i=(r+1)a_n-1$, then
    $Tf_i=n^i2^{(1-a_n/2)/\sqrt{a_n}}e_{(r+1)a_n}$. To express this in
    terms of the $f_i$'s, we notice that
    $f_{(r+1)a_n}=a_{n-r-1}(n^{(r+1)a_n}e_{(r+1)a_n}-e_0)$, which implies
    \marg{e:e(r+1)an}
    \begin{equation} \label{e:e(r+1)an}
      e_{(r+1)a_n}=n^{-(r+1)a_n}(a_{n-r-1}^{-1}f_{(r+1)a_n}+f_0),
    \end{equation}
    In this case $T^-f_i=0$.
    Analogously, if $i=a_n-1$, then
    $Tf_i=n^i2^{(1-a_n/2)/\sqrt{a_n}}e_{a_n}$. It follows from
    $f_{a_n}=a_{n-1}(n^{a_n}e_{a_n}-e_0)$ that
    \[
      Tf_i=n^{-1}2^{(1-a_n/2)/\sqrt{a_n}}(a_{n-1}^{-1}f_{a_n}+f_0),
    \]
    and again $T^-f_i=0$. Thus, case (B) produces no nontrivial
    entries in $T^-$.
  \item[C)] If $i$ is not the right end of the interval, i.e.
    $i<na_n+rb_n$, then $Tf_i=n^{-1}f_{i+1}$, so
    that $T^-f_i=0$. If $i=na_n+rb_n$, then
    \begin{multline*}
      Tf_i=
      n^{na_n+rb_n}e_{1+na_n+rb_n}-
      b_nn^{na_n+(r-1)b_n}e_{1+na_n+(r-1)b_n}\\
      =\varepsilon_1f_{1+na_n+rb_n}-\varepsilon_2f_{1+na_n+(r-1)b_n},
    \end{multline*}
    where $\varepsilon_1>0$ and
    $\varepsilon_2=b_nn^{-1}2^{(1+na_n-b_n/2)/\sqrt{b_n}}$. It follows that
    \marg{e:C}
    \begin{equation} \label{e:C}
      T^-f_{na_n+rb_n}=
      b_nn^{-1}2^{(1+na_n-b_n/2)/\sqrt{b_n}}f_{1+na_n+(r-1)b_n}.
    \end{equation}
  \item[D)] If $i<(r+1)(a_n+b_n)-1$, then
    $Tf_i=n^{-1}2^{1/\sqrt{b_n}}f_{i+1}$, so that $T^-f_i=0$. If
    $i=(r+1)(a_n+b_n)-1$ then
    \[
      Tf_i=n^i2^{(-a_n/2-(r+1)a_n/2+1)/\sqrt{b_n}}e_{(r+1)(a_n+b_n)}.
    \]
    Using (C) inductively we obtain the following identity:
    \begin{multline*}
      e_{(r+1)(a_n+b_n)}=n^{-(r+1)(a_n+b_n)}
      \{f_{(r+1)(a_n+b_n)}+b_nf_{(r+1)a_n+rb_n}+\dots\\
      +n_n^rf_{(r+1)a_n+b_n}\}
      +b_n^{r+1}n^{-(r+1)b_n}e_{(r+1)a_n}.
    \end{multline*}
    Substitute $e_{(r+1)a_n}$ from (\ref{e:e(r+1)an}) and notice that all the
    the coefficients are positive and, therefore, $T^-f_i=0$. Thus,
    case (D) does not produce any nontrivial entries in $T^-$.
\end{enumerate}

 Summarizing the calculations, the only nonzero entries of $T^-$ are
 given by (\ref{e:A}) and (\ref{e:C}):
  \[
    t^-_{1+v_{n-r},v_{n-r}+ra_n}=(n-r+1)^{-1-v_{n-r}}
        2^{(1+v_{n-r}-a_{n-r+1}/2)/\sqrt{a_{n-r+1}}}
        a_{n-r}(n-r)^{v_{n-r}}
  \]
  and
  \[
    t^-_{na_n+(r-1)b_n+1,na_n+rb_n}=b_nn^{-1}2^{(1+na_n-b_n/2)/\sqrt{b_n}}
  \]
  for all $0<r\le n$. To show that $T^-$ is nuclear it suffices to show
  $\sum_{k=0}^\infty\norm{T^-_{(k)}}_\infty<+\infty$. Look at the rows
  of $T^-$ containing non-zero entries.
  Notice that
  \[
    t^-_{1+v_{n-r},v_{n-r}+ra_n}\le
    a_{n-r}2^{(1+v_{n-r}-a_{n-r+1}/2)/\sqrt{a_{n-r+1}}}\le
    2^{-(1+v_{n-r})}
  \]
  for all $0<r\le n$ provided {\bf d} increases sufficiently
  rapidly. It follows that
  $\norm{T^-_{(1+v_m)}}_\infty\le 2^{-(1+v_m)}$ for every $m\ge 0$ and
  $\sum\limits_{m=0}^\infty \norm{T^-_{(1+v_m)}}_\infty\le
  \sum\limits_{m=0}^\infty 2^{-1-v_m}<1$.

  Further, the entries $t^-_{na_n+(r-1)b_n+1,na_n+rb_n}$ do not depend
  on $r$, and their contribution to
  $\sum_{k=0}^\infty\norm{T^-_{(k)}}_\infty$ does not
  exceed the sum of all of them, which can be easily estimated:
  \[
    \sum\limits_{n=1}^\infty\sum\limits_{r=1}^n
    b_nn^{-1}2^{(1+na_n-b_n/2)/\sqrt{b_n}}
    \le\sum\limits_{n=1}^\infty b_n2^{(1+na_n-b_n/2)/\sqrt{b_n}}
    \le\sum\limits_{n=1}^\infty 2^{-n}= 1,
  \]
  because
  $b_n2^{(1+na_n-b_n/2)/\sqrt{b_n}}\le 2^{-n}$ for all $n\ge 1$
  provided {\bf d} increases sufficiently rapidly.
  Thus, $\nu(T^-)\le\sum_{k=0}^\infty\norm{T^-_{(k)}}_\infty<2$
  provided {\bf d} increases sufficiently rapidly.
\end{proof}

\begin{corollary} \lab{c:main}
  C.~J.~Read's operator $T$ satisfies the following properties, provided
  {\bf d} increases sufficiently rapidly:
  \begin{enumerate}
  \item \lab{i:eigen} $\abs{T}$, $T^+$, and $T^-$ have positive eigenvectors;
  \item \lab{i:no ideals} Neither $\abs{T}$ nor $T^+$ has an invariant ideal.
  \end{enumerate}
\end{corollary}

\begin{proof}
  It follows from Theorem~\ref{t:general} and
  Lemma~\ref{l:nuclear} that $\abs{T}$ has a positive eigenvector.
  It was noticed in the proof of
  Lemma~\ref{l:nuclear} that $T^-f_0=0$, so that $T^-$ also has a
  positive eigenvector.

  To prove (\ref{i:no ideals}), assume that $J$ is a closed ideal in
  $\ell_1$ invariant under $\abs{T}$ or $T^+$, and that $0\neq x\in J$,
  then $x_k\neq 0$ for some $k\ge 0$, so that $f_k\in J$. It follows
  from the proof of Lemma~\ref{l:nuclear} that both $\abs{T}f_i$ and
  $T^+f_i$ have nonzero $(i+1)$-th component, implying $f_{k+1}\in J$.
  Proceeding inductively, we see that $f_i\in J$ for all $i\ge k$.
  Further, the proof of Lemma~\ref{l:nuclear} also shows that
  $(\abs{T}f_i)_0\neq 0$ and $(T^+f_i)_0\neq 0$ for infinitely many
  $i$'s, so that $f_0\in J$. It follows that $f_i\in J$ for every
  $i\ge 0$, so that $J=\ell_1$. In fact, (\ref{i:no ideals}) is a
  manifestation of the fact that a positive operator $S$ on
  $\ell_p$ $(1\le p<\infty)$ has no invariant ideals if and only if
  there is a path between every two columns of $S$
  (c.f.~\cite{AAB3,TV2}).

  It follows from~(\ref{i:no ideals}) and Theorem~\ref{t:aab} that $T^+$
  cannot be quasinilpotent. On the other hand, since $T^+=T+T^-$ then,
  analogously to the proof of Theorem~\ref{t:general}, we have
  $r_e(T^+)=0$. Then by Theorem~\ref{t:kr} we conclude that $r(T^+)$
  is a positive eigenvalue of $T^+$, corresponding to a positive
  eigenvector.
\end{proof}

  The last statement of Corollary~\ref{c:main} emphasizes that the
  hypothesis of not having invariant ideals in Theorem~\ref{t:general}
  is weaker than not having invariant subspaces. We do not know 
  if the analogues of the results of this paper hold
  for the operators produced in~\cite{R1,R2}.
  
  I would like to thank Prof.~Y.~A.~Abramovich for suggesting this
  problem to me and for our many discussions. I am thankful to
  Professors C.~D.~Aliprantis, V.~J.~Lomonosov, and C.~J.~Read for
  their interest in this work.

  We remark in conclusion that in~\cite{TV1}
  we use one of C.~J.~Read's operators to solve one more problem related
  to invariant subspaces. Namely, we construct operators $S_1$,
  $S_2$, and $K$ (not multiples of the identity) on $\ell_1$ such that $T$
  commutes with $S_1$, $S_1$ commutes with $S_2$, $S_2$ commutes with
  $K$, and $K$ is compact. This shows that the celebrated Lomonosov
  theorem cannot be extended to chains of four operators.

\end{document}